\documentclass[a4paper, 10pt, conference]{ieeeconf}      

\IEEEoverridecommandlockouts                              
\usepackage[utf8]{inputenc}
\usepackage[T1]{fontenc}
\usepackage{graphicx}      
\usepackage{amsmath}
\usepackage{amssymb}
\usepackage{esint}
\usepackage{multirow}
\usepackage{booktabs}
\usepackage{nicematrix}
\usepackage{tikz}
\usepackage{textcomp}
\usepackage{lipsum}
\makeatletter               
\let\NAT@parse\undefined    
\makeatother   
\usepackage{hyperref}   

\newcommand{\specialcell}[2][c]{%
  \begin{tabular}[#1]{@{}c@{}}#2\end{tabular}}
\newtheorem{assumption}{Assumption}
\newtheorem{definition}{Definition}
\newtheorem{remark}{Remark}
\DeclareMathOperator*{\argmin}{argmin}

\title{\LARGE \bf
Risk in Stochastic and Robust Model Predictive Path-Following Control for Vehicular Motion Planning
}

\author{Leon Tolksdorf$^{1, 2}$, Arturo Tejada$^{1, 3}$, Nathan van de Wouw$^{1}$ and Christian Birkner$^{2}$ 
\thanks{$^{1}$Department of Dynamics and Control, Eindhoven University of Technology, Eindhoven, The Netherlands, e-mail:
        {\tt\small \{l.t.tolksdorf, a.tejada.ruiz, n.v.d.wouw\}@tue.nl}}%
\thanks{$^{2}$CARISSMA Institute of Safety in Future Mobility, Technische Hochschule Ingolstadt, Ingolstadt, Germany, e-mail:
        {\tt\small \{leon.tolksdorf, christian.birkner\}@thi.de}}%
\thanks{$^{3}$TNO, Integrated Vehicle Safety, Helmond, The Netherlands, e-mail:
        {\tt\small arturo.tejadaruiz@tno.nl}}%
}

\newcommand\copyrighttext{%
  \footnotesize \textcopyright 2023 IEEE.  Personal use of this material is permitted.  Permission from IEEE must be obtained for all other uses, in any current or future media, including reprinting/republishing this material for advertising or promotional purposes, creating new collective works, for resale or redistribution to servers or lists, or reuse of any copyrighted component of this work in other works.}
\newcommand\copyrightnotice{%
\begin{tikzpicture}[remember picture,overlay]
\node[anchor=south,yshift=40pt] at (current page.south) {\fbox{\parbox{\dimexpr\textwidth-\fboxsep-\fboxrule\relax}{\copyrighttext}}};
\end{tikzpicture}%
}
\begin{document}

\maketitle

\thispagestyle{empty}
\pagestyle{empty}
\copyrightnotice
\begin{abstract}                
In automated driving, risk describes potential harm to passengers of an autonomous vehicle (AV) and other road users. Recent studies suggest that human-like driving behavior emerges from embedding risk in AV motion planning algorithms. Additionally, providing evidence that risk is minimized during the AV operation is essential to vehicle safety certification. However, there has yet to be a consensus on how to define and operationalize risk in motion planning or how to bound or minimize it during operation. In this paper, we define a stochastic risk measure and introduce it as a constraint into both robust and stochastic nonlinear model predictive path-following controllers (RMPC and SMPC respectively). We compare the vehicle's behavior arising from employing SMPC and RMPC with respect to safety and path-following performance. Further, the implementation of an automated driving example is provided, showcasing the effects of different risk tolerances and uncertainty growths in predictions of other road users for both cases. We find that the RMPC is significantly more conservative than the SMPC, while also displaying greater following errors towards references. Further, the RMPCs behavior cannot be considered as human-like. Moreover, unlike SMPC, the RMPC cannot account for different risk tolerances. The RMPC generates undesired driving behavior for even moderate uncertainties, which are handled better by the SMPC.
\end{abstract}

\begin{keywords}
autonomous vehicles, motion planning, path following, robust model predictive control, stochastic model predictive control, risk assessment
\end{keywords}

\section{Introduction}
Introducing autonomous vehicles (AVs) into traffic at scale will take a long period during which AVs and human-controlled vehicles will share the roads. This development leads to scenarios where AVs and human-controlled vehicles have to predict each other's future motion and interact \cite{Koopman.2017}. For AVs to be accepted by human drivers, they should display the kind of behaviors that human drivers expect from each other. Such behavior, once understood, could be incorporated into AVs and operationalized for motion planning.\\
Modeling human driving relies either on extensive use of artificial intelligence or on models for each individual aspect of driving \cite{Kolekar.2021}. Both approaches have limitations. Artificial intelligence lacks explainability and causal reasoning \cite{Marcus.02.01.2018}. Modeling each individual aspect of driving separately  leads to a fragmented motion planning design. However, in the quest to find an underlying theory of human driving, recent literature suggests that AVs can mimic human driving by maintaining a stochastic risk estimate below a threshold level. Further, a correlation between an objective risk estimate and perceived risk has been established \cite{Kolekar.2021}, indicating the feasibility of objective risk estimation methods for motion planning.

Here, in accordance with related literature, risk is understood to be a stochastic quantity that describes the potential harm resulting from driving with limited knowledge of current and future driving states. Using a precise risk definition, we develop a motion planning strategy that generates vehicle behavior by limiting risk along the motion plan\footnote{Current safety standards \cite{ISO26262,ISO21448} require car makers to show evidence that their vehicles operate with "acceptably" low risk. We believe the work presented here can support providing such evidence.}.

For motion planning, competing interests such as risk and travel time minimization must be balanced. Several motion planning methods are established (see, e.g., \cite{Gonzalez.2016}). We consider model predictive control (MPC) because it allows for the flexible integration of risk either in the objective function or as a constraint. We focus on path-following MPC formulations computing an actuator (e.g., velocity and turn rate) input sequence. Path-following is preferable over trajectory-tracking, because a time-independent reference path and a reference velocity can be computed offline first and the controller assigns the timing along the references online later \cite{Faulwasser.2016}. A path-following MPC scheme can effectively combine motion planning and control, balancing path-following performance (i.e., minimizing the error towards references) with risk minimization (potentially pushing the vehicle away from the references, i.e., around dynamic objects).

Two MPC formulations can accommodate the stochastic nature of risk: robust MPC (RMPC) or stochastic MPC (SMPC) \cite{Mayne.2016}. However, to the best of our understanding, the incorporation of risk in RMPC and SMPC schemes is still incipient. For instance, \cite{Batkovic.2022b} shows how to introduce stochastic constraints regarding uncertain positions of other road users into an extended trajectory-tracking RMPC formulation. It also provides conditions for recursive feasibility (which, in turn, guarantee AV safety) that rely on having uncertainties with bounded support. This allows the RMPC formulation to optimize for the worst-case (uncertainty) scenario, effectively sidestepping the stochastic nature of the uncertainties, potentially leading to more conservative behaviors \cite{Mayne.2014} (here, "conservative" means to take unnecessary large safety distances at the cost of path-following performance). Safety, in terms of recursive feasiblity for an SMPC is reported in \cite{Brudigam.2022} for a trajectory-tracking problem considering a linear system model. Safety is proven by assuming a safe backup planner that relies on worst-case scenario assumptions.\\
Note that unlike deterministic MPC, for which extensive tools exist to analyze reference convergence and recursive feasibility of trajectory-tracking and path-following problems \cite{Faulwasser.2012}, tools for analyzing RMPC and SMPC schemes for motion planning and control are limited. To the best of our knowledge, the integration of stochastic risk into nonlinear path-following SMPC and RMPC formulations is not reported in the literature, and recursive feasibility conditions for such formulations are only available for situations where risk can be upper-bounded by analyzing the worst-case scenarios. Furthermore, although the use of risk in motion planning has been reported (see, e.g., \cite{FlorianDamerow.2018,Hruschka.01.09.2021,Nyberg.7112021}), there is no consensus yet on its definition and operationalization. \\
This paper addresses some of these challenges. It provides a method to integrate stochastic risk as part of an MPC path-following problem formulation for automated driving applications. Herein, a general stochastic risk definition based on risk-inducing events is proposed. In addition, we compare two different approaches, SMPC and RMPC though a simulation-based case study of a representative automated driving scenario. The comparison is based on 1) safety, measured as the distance to an encountered object given a specific risk tolerance, and 2) path-following performance defined as reference error minimization. The simulation scenario consists of an AV that must follow a prescribed path in the presence of a moving obstacle. The effect of different rates of growing uncertainties in the predictions of the obstacle's motion and of various risk tolerances is also explored within the scenario.\\
The remainder of this article is organized as follows. Section \ref{riskdef} provides our risk definition. Section \ref{problem} shows how to integrate risk in the formulation of the path-following problem. Section \ref{rmpc} and \ref{smpc} derive the respective RMPC and SMPC formulations. In Section \ref{implementation} we present the simulation scenario, highlighting the required modelling elements. Section \ref{simulation} simulates the example with varying parameters. The results, with attention to our objectives, are discussed in \ref{discussion}. Section \ref{conclusions} highlights the benefits of SMPC over RMPC and presents an outlook for future research.

\textit{Notation:}
In our scenario, the AV (or ego vehicle) and the object are characterized by a configuration $\boldsymbol{y} := (\boldsymbol{q}, \theta) \in \mathcal{C}$. The configuration space $\mathcal{C}$ denotes the set of acceptable actor positions $\boldsymbol{q} := (c_1, c_2) \in \mathbb{R}^2$ and heading angles $\theta \in [0, 2\pi)$. 
A set of integers $\{a, a + 1, ..., b \}$, with $ a < b$, is denoted by $\mathbb{Z}_a^b$ and the set of positive reals including zero by $\mathbb{R}_{0,+}$. The identity matrix of size $p\times p$ is denoted by $\mathbb{I}_p$. We denote a configuration at time $k$ as $\boldsymbol{y}_{k}$ and a predicted configuration at time $n$ given information available at time $k$ as $\boldsymbol{y}_{n|k}$. The same notation is applied for states, inputs, and constraints. Finally, all variables associated with the ego vehicle and the object will be identified, respectively, with $e$ and  $o$ subscripts.

\section{Risk Definition}\label{riskdef}

The definitions of risk vary in the literature; nevertheless, safety standards \cite{ISO26262, ISO21448} formulate it as a combination of the probability of occurrence of harm and the severity of that harm\footnote{In a broader sense, the breaking of traffic laws, adhering to safety specifications \cite{Nyberg.7112021} and fairness among traffics participants \cite{Geisslinger.2021} could also be identified as risk-inducing events. However, this cannot necessarily be related to physical harm, shifting away from the notions of safety standards. To include risk-sources beyond safety standards, the definition of risk must be reformulated on the basis of a more general cost instead of physical harm. This will be addressed by future works of the authors.}. Examples for risk in the scope of these safety standards can be found in \cite{Tejada.2019, Vaicenavicius.}. Inspired by this approach, we define risk based on harm-inducing events $\mathcal{E}$, e.g., collisions or strong deceleration. However, our definition can be extended to other risk-inducing events.\\
Within a specific driving scenario, the occurrence of a harm-inducing event at time $k$ can be identified by ascertaining whether specific conditions among kinematic variables associated with the participants in the scenario are fulfilled at that time. These conditions can, for example, evaluate whether the physical boundaries of different participants overlap, representing a collision. Let $\boldsymbol{z}_{k} \in \mathbb{R}^{n_\mathcal{E}}$ denote the kinematic variables at time $k$ required to determine the occurrence of an event $\mathcal{E}$. Denote by $\mathcal{B}_{\mathcal{E},k}$ the subset in $\mathbb{R}^{n_\mathcal{E}}$ that satisfies the conditions used to identify the occurrence of $\mathcal{E}$ at time $k$. In general, $\boldsymbol{z}_{k}$ is a random vector with associated probability density function $p_{\boldsymbol{z},k}$, due to uncertainties involved in measuring or estimating the kinematic variables of interest. Lastly, let $s: \mathcal{B}_{\mathcal{E},k} \rightarrow \mathbb{R}_{+,0}$ denote a function that assigns the severity to every element of $\mathcal{B}_{\mathcal{E},k}$. Then, we can define the risk $R_k$ at time $k$ as
\begin{equation}\label{r1}
R_k := \mathbb{E}[s(\boldsymbol{z}_k)] = \int_{\mathcal{B}_{\mathcal{E},k}} s(\boldsymbol{z})p_{\boldsymbol{z},k}(\boldsymbol{z})\mathrm{d}\boldsymbol{z},
\end{equation}
that is, $R_k$ is interpreted as the expected severity of event $\mathcal{E}$ at time $k$. To constrain the risk at every time step, we introduce the risk tolerance $\epsilon \in \mathbb{R}_+$, such that $R_k \leq \epsilon$ represents an explicit risk constraint within our problem formulation. 

\section{Problem Formulation}\label{problem}
We proceed to integrate the risk defined in Section \ref{riskdef}, in the formulation of a path-following problem. Also, the control objectives and the control problem will be presented. While the problem setup introduces some simplification, we stress that our generic approach is not limited to the presented setup. 
\subsection{Problem Setup}\label{setup}
\begin{figure}
\begin{center}
\includegraphics[width=8.4cm]{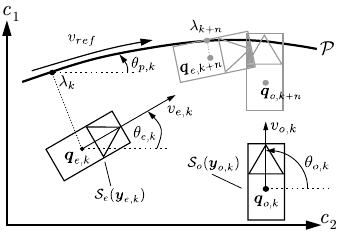}    
\end{center}
\caption{Illustration of the problem setup for two different time instances. At time $k$ the actors are displayed in black and for some time later $k+n$ in gray. Minimizing the reference error given an evolution of $\lambda$ on the path $\mathcal{P}$ would lead to a collision, where the gray area denotes $\mathcal{S}_e(\boldsymbol{y_{e,n+k}}) \cap \mathcal{S}_o(\boldsymbol{y_{o,n+k}})$.}
\label{problem setup}
\end{figure}

Consider a scenario consisting only of an AV, or ego vehicle, and another vehicle, called the object. The ego is provided with a reference path $\mathcal{P}$ and reference velocity $v_{ref} \in \mathbb{R}_{+, 0}$, and it is tasked with planning its own motion online to minimize the error with respect to both references within some finite horizon $\mathbb{Z}_{k_0}^{k_0 +K}$, with $K$ a non-negative integer. A path is defined as follows.
\begin{definition}\label{def1}(Path)
A path $\mathcal{P}$ is the image of a function $y_P: \mathbb{R} \rightarrow \mathcal{C}$ in the configuration space: $\mathcal{P} := \{\boldsymbol{y}_P(\lambda) \in \mathcal{C} \mid \lambda \in [\lambda_0, \lambda_g] \mapsto y_P(\lambda) \}$, where $y_P$ is sufficiently often differentiable with respect to $\lambda$.
\end{definition}
We assume that the path provided to the ego vehicle is a regular curve (i.e., its derivative never vanishes), satisfying Definition \ref{def1}. The path is exactly followable, that is, the path satisfies the physical limitations of the ego vehicle. Note that, we do not impose any collision constraints on the references. This means that exact \textit{following} of the path and velocity reference could lead the ego to collide with the object.
Here, a collision between the ego and the object is the only considered harm-inducing event $\mathcal{E}$. To check whether a collision occurred, we consider the ego vehicle's and object's configuration, respectively, denoted by $\boldsymbol{y}_{e, k} \in \mathcal{C}_{e,k} \subset \mathcal{C}$ and $\boldsymbol{y}_{o, k} \in \mathcal{C}_{o,k}\subset \mathcal{C}$ at time $k$. To estimate a collision's severity, we also require the respective scalar velocities $v_{e,k}, v_{o,k}$ pointing in the direction of the heading angles which are contained within the configurations. The velocities are elements of the respective time-dependent compact sets $\mathcal{V}_{e,k}, \mathcal{V}_{o,k} \subset \mathbb{R}$. Therefore, $\boldsymbol{z}_k = (\boldsymbol{y}_{e, k}^T, v_{e,k}, \boldsymbol{y}_{o, k}^T, v_{o,k})^T$ denote the kinematic variables required to determine the risk associated with a collision. 
A collision occurs if the occupied compact regions of both actors intersect. We denote the respective compact regions $\mathcal{S}_e(\boldsymbol{y}_{e, k}), \mathcal{S}_o(\boldsymbol{y}_{o, k}) \subset \mathbb{R}^2$  at time $k$. Hence, a collision is conditioned on $\mathcal{S}_e(\boldsymbol{y}_{e, k}) \cap \mathcal{S}_o(\boldsymbol{y}_{o, k}) \neq \emptyset$. Clearly, the collision condition depends on the configuration only, while the severity additionally requires the velocities. Thus we define the set $\mathcal{B}_{\mathcal{E}, k} := \{ \boldsymbol{z} \in \mathcal{C}_{e,k}  \times \mathcal{V}_{e,k} \times \mathcal{C}_{o,k} \times \mathcal{V}_{o,k}  \mid  \mathcal{S}_e(\boldsymbol{y}_{e}) \cap \mathcal{S}_o(\boldsymbol{y}_{o}) \neq \emptyset \}$ to be used in the risk definition in (\ref{r1}). The problem setup is illustrated in Figure \ref{problem setup}.

We make the following assumptions to clarify which information is available to the controller.

\begin{assumption}\label{ass1}(Object information) 
The object configuration $\boldsymbol{y}_{o,k}$ and velocity $v_{o, k}$ at time $k$ are measured (either by direct measurement or estimation) together with their associated independent probability densities $p_{\boldsymbol{y}_o, k}$ and $p_{v_o, k}$. For all future instances $n$ within a prediction horizon $\mathbb{Z}_k^{k+N}$, there is a prediction model available for the object, that propagates $(\boldsymbol{y}_{o, n|k}, v_{o, n|k})$ in the time-dependent set $\mathcal{C}_{o,n|k} \times \mathcal{V}_{o,n|k}$ and estimates the independent probability density functions $p_{\boldsymbol{y}_o, n|k}$ and $p_{v_o, n|k}$. 
\end{assumption}
\begin{assumption}\label{ass11}(Ego information) 
The kinematic variables associated with the ego vehicle are assumed to be measurable at all times without any uncertainty. 
\end{assumption}
We adopt Assumption \ref{ass1} since we are interested in understanding the influence of prediction uncertainties, captured within the risk measure, on the ego's behavior. The uncertainties associated with the object are generally more significant than those with the ego vehicle, justifying Assumption \ref{ass11}. 
\subsection{Control Objectives}
The problem is set up such that we seek to minimize an error towards time-independent references. However, the exact reference following may lead to collisions. To avoid collisions, we introduce a risk constraint that potentially forces the controller to deviate from the references, requiring the controller to balance reference following with collision risk. Thus, the controller performs online motion planning and vehicle actuation for the ego vehicle based on the following objectives:
\begin{itemize}
\item[(i)] minimize the following error to the reference path and velocity,
\item[(ii)] constraint satisfaction, explicitly including a risk constraint to enforce safety.
\end{itemize} 

\subsection{Control Problem}
In contrast to trajectory-tracking, path-following treats the time evolution along the path as an additional degree of freedom. Thereby path-following allows the controller to plan maneuvers by deviating from the reference velocity and path. Suppose that the motion dynamics of the ego are modeled by a discrete-time non-linear system of the form 
\begin{equation}\label{sys}
\begin{split}
\boldsymbol{x}_{k+1} &= f(\boldsymbol{x}_k, \boldsymbol{u}_{1,k}), \\
(\boldsymbol{y}_{e, k}^T, v_{e,k}) & = h(\boldsymbol{x}_{k}),
\end{split}
\end{equation}
where $\boldsymbol{x}_k \in \mathbb{R}^{n_x}$ and $\boldsymbol{u}_{1,k} \in \mathbb{R}^{n_{u_1}}$ represent the state and input vectors at time $k$ respectively. The function $f$ comprises the model equations and $h$ the output functions. The system is subject to state and input constraints, such that $\boldsymbol{x} \in \mathcal{X}$ and $\boldsymbol{u}_1 \in \mathcal{U}$ which must be satisfied at all times. The constraints represent actuator limitations and dynamic limits of the vehicle. We require the system to minimize the error with respect to the path, where path-following allows us to control the evolution along the path, parameterized by $\lambda \in [\lambda_0, \lambda_g]$ \cite{Faulwasser.2016}. To do so, an additional input  $u_{2,k} \in \mathcal{V}_{e,k}$ is introduced to control the evolution of $\lambda$ through the first-order difference equation
\begin{equation}\label{input}
\lambda_{k+1} = g(\lambda_k,\boldsymbol{y}_{e,k},u_{2,k}),
\end{equation}
where $g$ is a continuous function.

Next, we define the following error as 
\begin{equation}\label{error}
\begin{split}
\boldsymbol{e}_k & =
\begin{pmatrix}
h^T(\boldsymbol{x}_k) - \boldsymbol{y}_P^T(\lambda_k)\\
u_{2,k} - v_{ref}
\end{pmatrix},
\end{split}
\end{equation}
where $v_{ref} \in \mathbb{R}_{+,0}$ is the reference velocity.  Indeed, when $\boldsymbol{e}_k=\boldsymbol{0}$, the AV follows the reference path with the reference velocity. Note that we assume that only the objects configuration and velocity are uncertain, instead of all variables leading us to separate the kinematic variables as $\boldsymbol{z}_{k} = (\boldsymbol{z}_{e,k}, \boldsymbol{z}_{o})$.
Thus, $\boldsymbol{z}_{e,k} = (\boldsymbol{y}_{e,k}^T, v_{e,k})$ are deterministic variables\footnote{ $\boldsymbol{z}_{e,k}$ is a variable of the domain of integration, such that $\mathcal{B}_{\mathcal{E},k}(\boldsymbol{z}_{e,k}) := \{ (\boldsymbol{z}_{e,k}, \boldsymbol{z}_{o}), \boldsymbol{z}_{o} \in \times \mathcal{C}_{o,k} \times \mathcal{V}_{o,k}  \mid  \mathcal{S}_e(\boldsymbol{y}_{e,k}) \cap \mathcal{S}_o(\boldsymbol{y}_{o}) \neq \emptyset \}$.} and $\boldsymbol{z}_{o} = (\boldsymbol{y}_{o}^T, v_{o})$ are random variables. Further, given $s: \mathbb{R}^8 \rightarrow \mathbb{R}_{0,+}$, the general risk from (\ref{r1}) becomes
\begin{equation}\label{rr1}
R_k = \iiint\displaylimits_{\mathcal{B}_{\mathcal{E},k}(\boldsymbol{z}_{e,k})}s(\boldsymbol{z}_k)\delta_{\boldsymbol{z}_{e,k}}(\boldsymbol{z}_e) p_{\boldsymbol{y}_o,k}(\boldsymbol{y})p_{v_o,k}(v)\text{d}\boldsymbol{z}_e\text{d}\boldsymbol{y}\text{d}v,
\end{equation}
where $p_{\boldsymbol{z},k} = \delta_{\boldsymbol{z}_{e,k}} p_{\boldsymbol{y}_o,k}p_{v_o,k}$ due to the independence assumption, and $\delta_{\boldsymbol{z}_{e,k}}$ is the Dirac distribution over $\boldsymbol{z}_{e,k}$.
\textit{Problem:} Given the risk measure (\ref{rr1}), the system dynamics (\ref{sys}) - (\ref{input}), a path $y_P$ satisfying Definition \ref{def1} and a reference velocity $v_{ref}$, design a controller generating inputs $\boldsymbol{u}_1, u_2$ such that:

\begin{itemize}
\item[(a)] the error in (\ref{error}) is minimized, in the sense that $\sum_{k =k_0}^{k_0 + K} \| \boldsymbol{e}_k\|$ is minimized,
\item[(b)] the constraints $\boldsymbol{x}_k \in \mathcal{X}$, $(\boldsymbol{u}_{1,k}, u_{2,k}) \in \mathcal{U} \times \mathcal{V}_{e,k}$, and $R_k \leq \epsilon$ are satisfied for all $k$.
\end{itemize}
In the next two sections, we discuss two distinct MPC approaches to solving this problem.
\section{Robust Model Predictive Control}\label{rmpc}

We take the RMPC approach of \cite{Batkovic.2022b}, formulated for an extended trajectory-tracking problem with an uncertain constraint, as the basis for tackling the problem stated in Section \ref{problem}. 
In RMPC, the worst-case scenario is found by realizing the random variables such that cost or a constrained value are maximized. Given our problem formulation, the ego vehicle is deterministic; hence, the path-following error is deterministic, leaving the uncertainty only within the risk constraint, defined in the following.
\begin{definition}\label{def2}(Worst-case scenario: risk constraint)
The worst-case scenario given by (the kinematic variables) $\boldsymbol{z}\in \mathcal{B}_{\mathcal{E},n|k}(\boldsymbol{z}_{e,n|k})$ that maximizes the severity function $s$.
\end{definition}
As a consequence of Definition \ref{def2}, if one further assumes that $s$ is continuous and $\mathcal{C}_{o,n|k} \times \mathcal{V}_{o,n|k}$ is compact, the worst-case risk is given by 
\begin{equation}\label{wca}
R_{wc, n|k} = \max_{\boldsymbol{z}\in \mathcal{B}_{\mathcal{E},n|k}(\boldsymbol{z}_{e,n|k})} s(\boldsymbol{z}).
\end{equation}
Note that this only holds under Assumptions \ref{ass1} and \ref{ass11} and the assumption of compactness of $\mathcal{C}_{o,n|k} \times \mathcal{V}_{o,n|k}$, which guarantee that $\mathcal{B}_{\mathcal{E},n|k}(\boldsymbol{z}_{e,n|k})$ is also compact. The worst-case risk might be computeable in more general cases. 
Thus, a careful construction of $\mathcal{B}_{\mathcal{E},n|k}(\boldsymbol{z}_{e,n|k})$ and random variables with compact support are required to operationalize a general risk measure within RMPC. 
The cost function to be minimized at each sampling instance is
\begin{equation}
\begin{split}
J(\boldsymbol{y}_{e,k}, \lambda_k, \boldsymbol{u}_1, u_2)= E(\boldsymbol{y}_{e, k + N|k}, \lambda_{k+N|k})\\
 + \sum_{n = k}^{k + N-1}F(\boldsymbol{e}_{n|k}, \lambda_{n|k}, \boldsymbol{u}_{1, n|k}, u_{2,n|k})\label{r},
\end{split}
\end{equation}
where $\boldsymbol{u}_{1}, u_{2}$ denote input sequences over the prediction horizon of length $N$. Further, $E: \mathcal{C} \times [\lambda_0, \lambda_g] \rightarrow \mathbb{R}_{+,0}$ is the terminal cost and $F: \mathbb{R}^4 \times [\lambda_0, \lambda_g] \times \mathcal{U}\times \mathcal{V}_e \rightarrow \mathbb{R}_{+,0}$ represents the stage cost. Note that this is a generic formulation. In practice, the stage and terminal cost might assign a cost to fewer variables. To achieve our control objectives, it is sufficient to minimize the error. The path-following RMPC formulation is stated as
\begin{subequations}
\begin{equation}
\begin{split}
&V^{RMPC}(\boldsymbol{y}_{e,k}, \lambda_k) := \min_{\boldsymbol{u}_1, u_2} J(\boldsymbol{y}_{e,k}, \lambda_k, \boldsymbol{u}_1, u_2)\label{ra}, \\
\end{split}
\end{equation}
subject to: 
\begin{align}
 & \;\;\;\; \boldsymbol{x}_{k|k} = \boldsymbol{x}_k, \;\;\;\; \lambda_{k|k} = \lambda_k,\label{rb} \\
& \forall n \in \mathbb{Z}_k^{k+N-1}: \boldsymbol{x}_{n+1|k} = f(\boldsymbol{x}_{n|k}, \boldsymbol{u}_{1,n|k}),  \label{rc}\\ 
& \;\;\;\; (\boldsymbol{y}_{e,n|k}^T, v_{e, n|k})  = h(\boldsymbol{x}_{n|k}),\label{rd}\\
& \;\;\;\; \lambda_{n+1|k} = g(\lambda_{n|k},u_{2,{n|k}}), \label{re}\\ 
& \;\;\;\; \boldsymbol{e}_{n|k} = (h^T(\boldsymbol{x}_{n|k}) - \boldsymbol{y}_P^T(\lambda_{n|k}), u_{2,n|k} - v_{ref})^T,\label{rf} \\
&\;\;\;\;  \boldsymbol{x}_{n|k} \in \mathcal{X}, \;\;\;\;(\boldsymbol{u}_{1, n|k}, u_{2, n|k}) \in \mathcal{U} \times \mathcal{V}_e, \label{rg}\\ 
&  \forall n \in \mathbb{Z}_k^{k+N}: \lambda_n \in [\lambda_0, \lambda_g], \label{rh}\\
&\;\;\;\; R_{wc, n|k} \leq \epsilon,\label{ri}\\
& \;\;\;\; (\boldsymbol{x}_{k + N|k}, \lambda_{k+N|k}) \in \mathcal{T} \subset \mathcal{X} \times [\lambda_0, \lambda_g] \label{rj}.
\end{align}
\end{subequations}
To enforce that the prediction starts at the current states, the constraint (\ref{rb}) is implemented. The system model is integrated into (\ref{rc}) - (\ref{re}). The error is given in (\ref{rf}) to support path-following by minimization of this error. State, input, and risk constraints are implemented in (\ref{rg}) and (\ref{rh}). The worst-case risk constraint is provided in (\ref{ri}). Lastly, we we demand that the state and path variable reach a terminal region $\mathcal{T}$ in (\ref{rj}). Due to Assumption \ref{ass11} the objective function (\ref{r}) is not a worst-case cost (it just is a deterministic cost), and the system model (\ref{rc}) - (\ref{re}) remains without any uncertainty.

\begin{remark}\label{rem1}(Conditions for recursive feasibility and asymptotic stability for RMPC)
To prove recursive feasibility and asymptotic stability, conditions on the system, cost, constraints, and references are required. Also, stabilizing terminal conditions need to be assumed. Strictly speaking, the RMPC scheme (\ref{ra}) - (\ref{rj}) has not yet been proven to be recursively feasible or asymptotically stabilize. Still, by assessment of the authors, the proof of recursive feasibility in \cite{Batkovic.2022b} for an extended trajectory-tracking problem in conjunction with the MPC path-following work of \cite{Faulwasser.2012} can be adapted without much modification to the problem at hand.
\end{remark}

\section{Stochastic Model Predictive Control}\label{smpc}

A different approach to treating the underlying uncertainty in the risk definition is based on SMPC: Here, one must not only use the uncertainty realization associated to the worst-case scenario, which essentially convert the stochastic problem into a deterministic problem (as done in the RMPC approach); instead, the stochastic problem setting can be directly integrated into the SMPC formulation. 

For the SMPC formulation, we do not require the worst-case scenario (see Definition \ref{def2}). We employ (\ref{r}) as a cost function and keep the formulation (\ref{rb}) - (\ref{rj}) except for (\ref{ri}), where instead of the worst-case risk constraint $ R_{wc, n|k} \leq \epsilon$ we implement the risk Equation \ref{rr1} directly as a stochastic risk constraint $ R_{n|k} \leq \epsilon$. Due to Assumption \ref{ass11}, the objective function (\ref{r}) is not an expected value. Therefore we obtain an optimal control problem with a chance constraint of the expectation type, which reads
\begin{equation}
V^{SMPC}(\boldsymbol{y}_{e,k}, \lambda_k):= \min_{\boldsymbol{u}_1, u_2} J(\boldsymbol{y}_{e,k}, \lambda_k, \boldsymbol{u}_1, u_2),
\end{equation}
subject to (\ref{rb}) - (\ref{rh}), $R_{n|k} \leq \epsilon$, (\ref{rj}).
\begin{remark}\label{rem2}(Conditions for recursive feasibility and asymptotic stability for SMPC)
Similar to Remark \ref{rem1}, the SMPC scheme has not yet been proven to be recursively feasible nor asymptotically stabilize. Again,  conditions on the system, costs, constraints, references, as well as stabilizing terminal conditions are required. If recursive feasibility is of interest (following the definition of safety by \cite{Batkovic.2022b}), the problem must be reformulated to a linear system model so that existing proofs (see, e.g., \cite{Hewing.2020b, Schluter.2022, Kohler.2022}) for recursive feasibility can be adapted to the path-following problem at hand. However, this requires further assumptions on the uncertainty, and it is an open research challenge that will be addressed by future work of the authors. 
\end{remark}

\section{Risk-based Path-Following Control of an AV}\label{implementation}

Embedding the RMPC and SMPC into a simulation environment or real-world application requires extensive modeling. We choose to keep the individual models of low complexity, i.e., simplified dynamic model, circular shape approximations and linearly growing Gaussian uncertainty, such that the results are easily traceable and the impact of different risk tolerances and uncertainty growths on the AV's behavior is explicit. While our example works with low-complexity models, we stress that the generic approach will also apply to higher-fidelity models.

\subsection{System Model}

For this case study, a unicycle vehicle model is selected, which represents all configuration variables. This model ignores the correlation between the velocity and turn rate as an input and assumes moving along a curved trajectory. Since the state vector equals the configuration vector, the output equation $h$ in (\ref{sys}) is not required; the model reads:
\begin{equation*}\label{i1}
\dot{\boldsymbol{x}}(t) = 
\begin{pmatrix}
\dot{c}_1(t)\\
\dot{c}_2(t)\\
\dot{\theta}(t)
\end{pmatrix}=
\begin{pmatrix}
v_e(t)\cos{(\theta_e(t))}\\
v_e(t)\sin{(\theta_e(t))}\\
\omega_e(t)
\end{pmatrix},
\end{equation*}
where $t \in \mathbb{R}_{+,0}$ represents continuous time; the velocity $v_e \in \mathbb{R}$ and turn rate $\omega_e \in \mathbb{R}$ are the inputs to the model. We use a forward Euler discretization with the step size $T$, yielding the discrete-time model as 
\begin{equation}\label{i3}
\boldsymbol{x}_{k+1} = \boldsymbol{x}_{k} + 
\begin{pmatrix}
\frac{v_{e,k}}{\omega_{e,k}}(\sin{(\theta_{e,k}+\omega_{e,k} T)}-\sin{(\theta_{e,k})})\\
\frac{v_{e,k}}{\omega_{e,k}}(\cos{(\theta_{e,k})}-\cos{(\theta_{e,k}+\omega_{e,k} T)})\\
\omega_{e,k} T
\end{pmatrix}.
\end{equation}
\subsection{Path and Timing Law}
\begin{figure}
\begin{center}
\includegraphics[width=9.5cm]{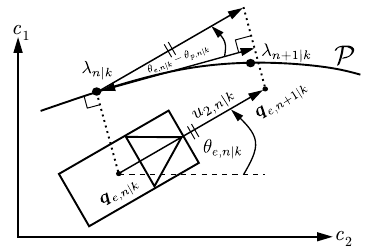}    
\end{center}
\caption{Approximation of path reference points $\lambda_{n|k}$. The dotted lines represent the projection of the displacement vector $u_{2,n|k}$ onto the nearest point on the path $\lambda_{n|k}$. The length of the tangent, given the relative angle $\theta_{e, n|k} - \theta_{p, n|k}$, is then utilized to calculate $\lambda_{n+1|k}$. One can see that estimated point $\lambda_{n+1|k}$ is only an approximation of the actual closest point from the ego vehicle to the path at $n+1|k$. 
}
\label{control law}
\end{figure}

As a path we utilize a regular curve, satisfying Definition \ref{def1} with constant curvature $\kappa$ along an interval $\lambda \in [\lambda_0, \lambda_g]$. For reference feasibility (i.e., the curve is followable), $\kappa$ must be less than the turn rate limit of the vehicle model. The initial point on the path is found by the closest point to the path from the initial vehicle configuration, that is
\begin{equation*}\label{path}
\lambda_{k|k} = \argmin_{\lambda \in [\lambda_0, \lambda_g]} \| \boldsymbol{y}_{e, k|k} - \boldsymbol{y}_P(\lambda) \|.
\end{equation*}
Finding the reference point by minimization is computationally intensive. Therefore, after computing the initial point $\lambda_{k|k}$, we approximate the following points within the prediction horizon by projecting the displacement vector, travelled by the ego vehicle in $T$ seconds due to the velocity $u_{2,n|k}$, onto the tangent to the path at $\lambda_{n|k}$. Thus, the next point is found with
\begin{equation*}
    \lambda_{n+1|k} = \lambda_{n|k} + u_{2,n|k}\cos{(\theta_{e, n|k} - \theta_{p, n|k})}T,
\end{equation*}
where $\theta_{p, n|k}$ represents the reference heading angle, i.e,. the angle of the tangent to the path at $\lambda_{n|k}$, see Figure \ref{control law}.

\subsection{Risk Model}
Given the problem setup, we are required to derive the risk of collision. We can assign a risk value to a collision using a shape approximation, severity model, and a description of the uncertainty. 
\subsubsection{Shape Approximation} The actor's shapes are assumed to be circles. Consequently, their orientation can be neglected while determining inter-actor distances. Thus, the actor's orientation within the configuration is not required in the following. Recall that $\boldsymbol{q}_k = (c_{1,k}, c_{2,k})^T$ is the geometric center of an actor at time instance $k$. The actor occupies a circular region in $\mathbb{R}^2$ of radius $r \in \mathbb{R}_+$. Hence 
\begin{equation}\label{shape}
\mathcal{S}(\boldsymbol{y}_k) = \{ (c_1, c_2)^T \in \mathbb{R}^2 \mid \| \boldsymbol{q}_k - (c_1, c_2)^T\| \leq r \}.
\end{equation} 
Approximating the actor shapes by circles reduces the checking for collision to examining if the Euclidean distance between the centers of both circles is less or equal then the sum of both radii; $ \|\boldsymbol{q}_{e} - \boldsymbol{q}_{o}\| \leq r_e + r_o$, where $r_e, r_o$ denote the respective radii of the ego vehicle and object. Reducing the dimension of of $\boldsymbol{z} \in \mathbb{R}^8$ by the respective heading angles is indicated by $\boldsymbol{\hat{z}}\in \mathbb{R}^6$. Therefore, we obtain $\mathcal{\hat{B}}_{\mathcal{E},k}(\boldsymbol{\hat{z}}_{e,k}) := \{ (\boldsymbol{\hat{z}}_{e,k}, \boldsymbol{\hat{z}}_{o}),\boldsymbol{z}_{o} \in \mathbb{R}^2 \times \mathcal{V}_{o, k} \mid   r_e + r_o \geq \|\boldsymbol{q}_{e, k} - \boldsymbol{q}_{o} \| \}$. 
\subsubsection{Severity Model}
The literature suggest various models to estimate the collision severity (see, e.g., \cite{Hruschka.01.09.2021}). Due to its computational efficiency, we use the differential kinetic energy as a severity model $s: \mathcal{V}_{e,k}\times \mathcal{V}_{o,k} \rightarrow \mathbb{R}_{0,+}$, reading
\begin{equation}\label{s}
s(v_{e,k}, v_{o,k}) = \frac{1}{2} \left| \left( m_e v_{e,k}^2 - m_o v_{o,k}^2 \right) \right|,
\end{equation}
where $m_e, m_o \in \mathbb{R}_+$ are the masses of the respective actors\footnote{Note that we reduced the domain of $s$ from $\mathcal{B}_{\mathcal{E},k}(\boldsymbol{z}_{e,k}) \subset \mathbb{R}^8$, see (\ref{rr1}), to $\mathbb{R}^2$ by not requiring the positions and orientations of the actors.}. 
\subsubsection{Uncertainty Model}
With (\ref{shape}) and (\ref{s}), the risk depends on the position and velocity. In the problem setup (see Section \ref{setup}), we assumed consistent configuration and velocity pairs with associated probability density functions for the object as given (see Assumption \ref{ass1}). Additionally, we assumed perfect knowledge of the ego vehicle's current and future kinematic variables (see Assumption \ref{ass11}). Furthermore, for the circular shape approximation of both actors as well as for the severity model we do not require the respective heading angles. Therefore the remaining random variables are the object's position $\boldsymbol{q}_o$ and velocity $v_o$. We choose truncated Gaussians centered around $\boldsymbol{q}_{o,k}$ and $v_{o,k}$ as the distributions with independent probability density functions $p_{\boldsymbol{q}_o, k}, p_{v, k}$ respectively at time $k$. The truncation bounds at time $k$ are $\boldsymbol{q}_{min, k} = (c_{1, min,k}, c_{2, min,k})^T, \boldsymbol{q}_{max, k}= (c_{1, max,k}, c_{2, max,k})^T$ for the object's position and $v_{min, k}, v_{max, k}$ for the object's velocity. $\sigma_{c_1,k}, \sigma_{c_2,k}, \sigma_{v, k}$ are the associated, respective standard deviations composed in $\Sigma_{k} =Diag(\sigma_{c_1,k}, \sigma_{c_2,k}, \sigma_{v, k})$. To model a growing uncertainty along the prediction horizon, we allow the standard deviation and truncation bounds to grow linearly over time by introducing an additive term for each as follows: 
\begin{align*}
&\Sigma_{n+1|k} = \Sigma_{n|k} + Q,\\
&\boldsymbol{q}_{min, n+1|k} = \boldsymbol{q}_{min, n|k} - \Delta \boldsymbol{q},\\
&\boldsymbol{q}_{max, n+1|k} = \boldsymbol{q}_{max, n|k} + \Delta \boldsymbol{q}, \\
&v_{min, n+1|k} = v_{min, n|k} - \Delta v,\\
&v_{max, n+1|k} = v_{max, n|k} + \Delta v,  
\end{align*} $ \forall n \in \mathbb{Z}_k^{k+N}$. The additive terms $Q \in \mathbb{R}^{3 \times 3}$, $\Delta \boldsymbol{q} \in \mathbb{R}^2_{0,_+}$, and $\Delta v \in \mathbb{R}_{0,_+}$ are constant. 

\begin{remark}
Note that we chose a truncated Gaussian so that all random variables are bounded. The bounding, however, is only required for the RMPC (see Definition \ref{def2}). Nevertheless, we apply the same uncertainty model to the SMPC and RMPC to ensure consistency within the comparison.
\end{remark}
\subsubsection{Risk Derivation}
With the characterized uncertainty, we introduce (\ref{shape}) and (\ref{s}) into (\ref{rr1}), we obtain $R_k =$

\begin{equation} \label{risk2}
 \iiint\displaylimits_{\mathcal{\hat{B}}_{\mathcal{E},k}(\boldsymbol{\hat{z}}_{e,k})}s(v_{e,k}, v_{o}) \delta_{\boldsymbol{\hat{z}}_{e,k}}(\boldsymbol{\hat{z}}_e) p_{v, k}(v_o) p_{\boldsymbol{q}_o, k}(\boldsymbol{q}_o)\text{d}\boldsymbol{\hat{z}}_e\text{d}\boldsymbol{q}_o\text{d}v_o. \\
\end{equation} 
To compute (\ref{risk2}), we apply Monte Carlo sampling (MCS). In order to numerically implement MCS, we utilize the collision indicator function
\begin{equation}\label{indicator}
I_C(\boldsymbol{q}_{e,k}, \boldsymbol{q}_{o,k}) = 
  \begin{cases}
    1       & \quad \text{if } \|\boldsymbol{q}_{e, k} - \boldsymbol{q}_{o, k}\| \leq  r_e + r_o, \\
    0  & \quad \text{otherwise.}
  \end{cases}
\end{equation}
With the indicator function (\ref{indicator}), the risk (\ref{risk2}) is approximated by the law of large numbers, where $(\boldsymbol{q}_{o,j,k}, v_{o, j,k})$ is one of $J$ samples drawn from the densities $p_{\boldsymbol{q}_o, k}, p_{v, k}$. This leads the risk (\ref{risk2}) to be approximated as
\begin{equation}\label{MCS}
\begin{split}
R_k \approx \frac{1}{J}\sum_{j = 1}^{J} I_C(\boldsymbol{q}_{e,k}, \boldsymbol{q}_{o,j,k})s(v_{e,k}, v_{o, j, k}).
\end{split}
\end{equation}
\subsection{Robust Risk Implementation}
For the robust MPC implementation, the challenge is to find the object's position and velocity, maximizing the risk. This yields the optimization problem
\begin{subequations}
\begin{align}\label{ria}
R_{wc, n|k} = \max_{\boldsymbol{q}_{o}, v_o}  I_C(\boldsymbol{q}_{e,n|k}, \boldsymbol{q}_{o})s(v_{e,n|k}, v_{o}), \text{s.b.t.:}\\
\boldsymbol{q}_{o} \in [ \boldsymbol{q}_{min, n|k},  \boldsymbol{q}_{max, n|k}], v_{o} \in [v_{min, n|k}, v_{max, n|k}]. \label{rib}
\end{align}
\end{subequations}
While (\ref{ria}) maximizes the severity and checks for a collision with the indicator function, constraints (\ref{rib}) enforce a solution under the given prediction (see Assumption \ref{ass1}). Due to the indicator function in (\ref{ria}), the constraint function is nonsmooth, ruling out gradient-based solvers. We propose to under-approximate the worst-case risk, arguing that the RMPC will be (even) more conservative with the actual worst-case. Hence, we find such under-approximation of the risk by a brute-force approach using a uniform sampling of $L$ points within each truncation interval, such that
\begin{equation} \label{ric}
\begin{split}
&(\boldsymbol{q}_o, v_o) \in [c_{1, min, n|k}, ..., c_{1, max, n|k}] \times \\
& [c_{2, min, n|k}, ..., c_{2, max, n|k}] 
\times [v_{min, n|k}, ..., v_{max, n|k}].
\end{split}
\end{equation}
We replace (\ref{ria}) by a maximization of the risk over the grid points in (\ref{ric}).

\subsection{Stochastic Risk Implementation}
For the SMPC implementation, we use random sampling to compute the risk from $(\ref{MCS})$. Hence, we directly implement (\ref{MCS}) within $V^{SMPC}(\boldsymbol{y}_{e,k}, \lambda_k)$ for $J$ samples.

\subsection{Cost Function and Terminal Conditions}
In the theoretical RMPC $V^{RMPC}(\boldsymbol{y}_{e,k}, \lambda_k)$ and SMPC $V^{SMPC}(\boldsymbol{y}_{e,k}, \lambda_k)$ formulation, terminal cost and constraints are proposed. In practice, these increase the complexity; the construction of terminal regions (\ref{ri}) has substantial implications on the recursive feasibility since one must construct a terminal region that can always be reached without constraint violation. The inclusion of the risk measure within the cost function is also possible, but similarly requires assumptions to ensure recursive feasibility. This is due to reference error minimization might increase the risk, and thus a monotonous decrease in cost while approaching the references is not provided, violating common assumptions. We leave this challenge for future work, where the recursive feasibility of the schemes will be investigated.
The cost for this work is a quadratic stage cost on the error to the references and we do not employ terminal cost, hence
\begin{equation*}\label{cf1}
F(\boldsymbol{e}_{n|k}) = \boldsymbol{e}_{n|k}^T W \boldsymbol{e}_{n|k},
\end{equation*}
where $W \in \mathbb{R}^{4 \times 4}$ is a positive definite weighting matrix. 

\begin{figure*}[t!]
\includegraphics{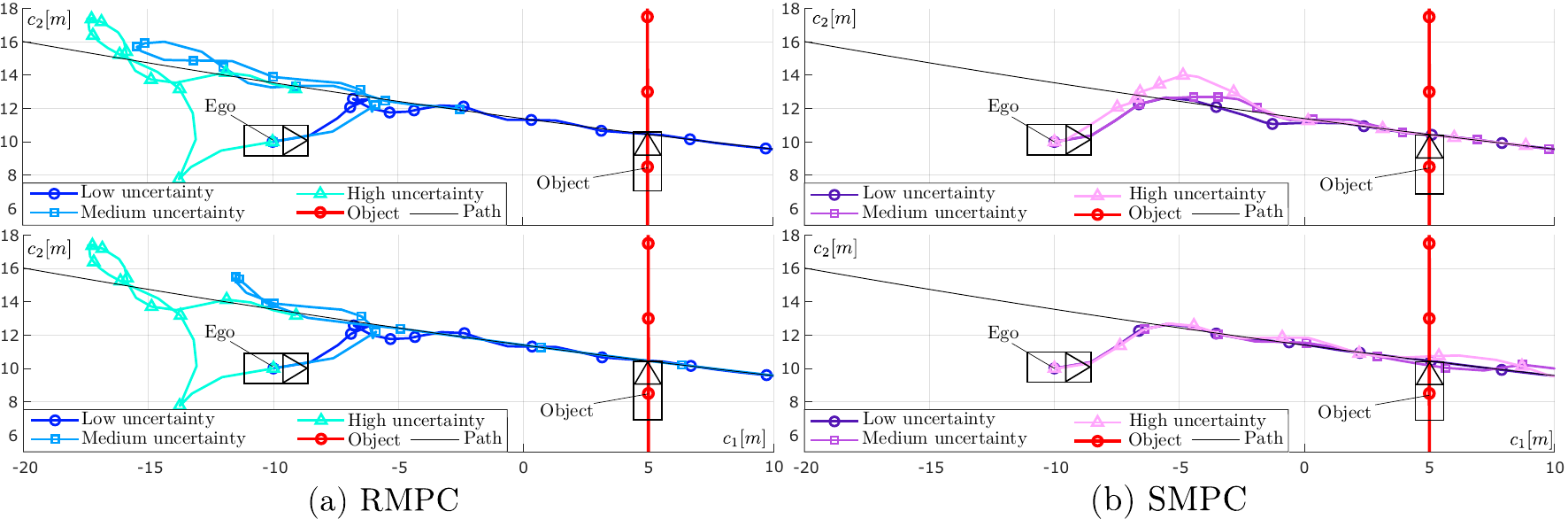}   
\caption{RMPC and SMPC comparison for three different uncertainties and two risk tolerances. Top plots: $\epsilon = 0$, bottom plots: $\epsilon = 2500$. In all plots the ego vehicle is displayed at $k=0s$ and the object at $k=9s$, the markers represent a time difference of $1.5 s$. The RMPC (a) takes significant deviations from the reference, increasing with the level of uncertainty for both risk tolerance settings. For all cases, the RMPC avoids a collision. The SMPC (b) requires less time to converge to the path and leaves less distance to the object while, in all cases, successfully avoiding a collision. Further, the SMPC leaves less distance to the object when the risk threshold increases. Also, under increasing uncertainty, the SMPC leaves more distance to the object.}
\label{comp} 
\end{figure*}
\section{Simulation Case Study}\label{simulation}

To demonstrate the impact of risk on the RMPC and SMPC formulations, we design a scenario where the ego is initially not on the reference path, and the initial ego velocity is set to zero. At some distance, one object crosses the reference path such that it would collide with the ego if the ego vehicle would closely follow its initial reference (comp. Figure \ref{problem setup}). We simulate this scenario with varying risk thresholds $\epsilon$. Additionally, we simulate different growth rates $Q, \Delta \boldsymbol{q} = (\Delta c_1, \Delta c_2), \Delta v$. While the complete list of parameters is provided in the Appendix, the varying uncertainty parameters are displayed in Table \ref{parameters}. For each of the uncertainty settings, that is; low, medium, and high, all six risk tolerances (see Table \ref{results}) are simulated with both methods. In total, $36$ simulations are conducted. 

\begin{figure}
\begin{center}
\includegraphics[width=8.4cm]{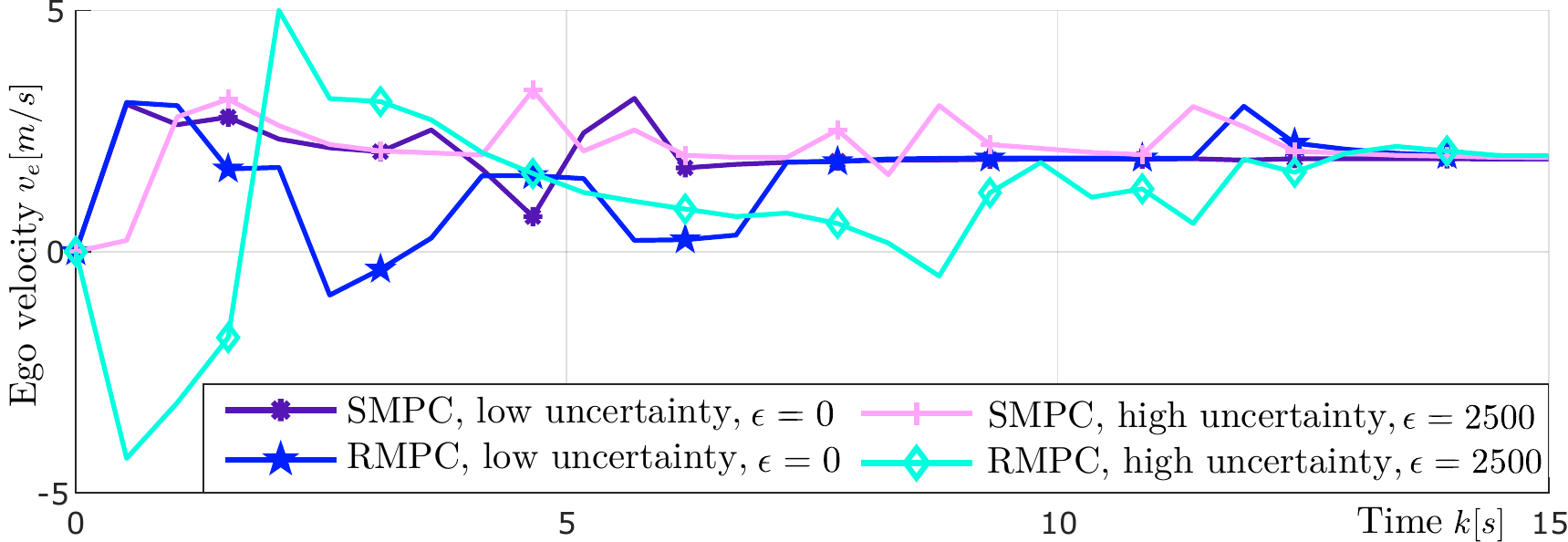}    
\end{center}
\caption{Respective ego vehicle velocities for some of the cases presented in Fig. \ref{comp}. The SMPC converges earlier than the RMPC to the reference velocity of $v_{ref} = 4 $ in both settings. Under high uncertainty, the RMPC chooses negative velocities, while the SMPC keeps positive velocities.}
\label{velo}
\end{figure}

\begin{table}[hb]
\begin{center}
\caption{Varying Simulation Parameters}\label{parameters}
\begin{tabular}{cccc}
Parameters & Low & Medium & High \\\hline
Additive uncertainty $Q$ &  $0.1 \cdot \mathbb{I}_3$ & $0.8\cdot \mathbb{I}_3$ & $1.5\cdot \mathbb{I}_3$ \\ 
\specialcell{Truncation growth\\$(\Delta c_1, \Delta c_2, \Delta v) $ }& $(1 ,1,  1) $ & $( 2, 2, 2)  $&$(3 , 3 , 3) $\\ \hline
\end{tabular}
\end{center}
\end{table}
Figure \ref{comp} presents the simulation for the lowest and highest risk tolerance. Figure \ref{velo} displays the ego velocities generated by the respective controllers for selected cases presented in Figure \ref{comp}. To compare both approaches, we measure the path-following performance and safety. The former is measured as the accumulated error over the simulation time, i.e.  $e_{acc} = \sum_{k =k_0}^{k_0 + K} \| \boldsymbol{e}_k\|$, whereas safety is measured as the minimum Euclidean distance $d_{min}$ from the ego's center to the object's center, which if $d_{min} \leq r_e + r_o$ indicates a collision. The results for all simulated scenarios are presented in Table \ref{results}. Note that $d_{min} = 3.0$ denotes a rounded value of $d_{min} > 3$.

\begin{table}[hb]
\begin{center}
\caption{Results}\label{results}
\begin{NiceTabular}{cc|cc|cc|cc}
 & &\multicolumn{2}{c}{Low uncert.}&\multicolumn{2}{c}{Medium uncert.} & \multicolumn{2}{c}{High uncert.}\\ 
 \multicolumn{2}{c}{} &$e_{acc}$& $d_{min}$& $e_{acc}$& $d_{min}$&$e_{acc}$& $d_{min}$ \\ \hline
 \multicolumn{1}{l}{\multirow{6}{*}{\rotatebox[origin=c]{90}{RMPC}}} &
   $\epsilon = 0$       &53.5 &9.1  &72.3 &14.7 &98.7 &20.3 \\
&$\epsilon = 500$      &53.5 &9.1  &72.2 &14.7 &98.7 &20.3 \\ 
&$\epsilon = 1000$    &53.5 &9.1   &72.5 &14.7 &98.7 &20.3\\ 
&$\epsilon = 1500$    &53.5 &9.1  &71.7  &14.7   &98.7 &20.3\\
&$\epsilon = 2000$    &53.5 &9.1  &70.9 &14.7  &98.8 &20.3\\
&$\epsilon = 2500$    &53.5  &9.1  &71.2 &14.7   &98.7 &20.3\\  \hline
 \multicolumn{1}{l}{\multirow{6}{*}{\rotatebox[origin=c]{90}{SMPC}}} & 
  $\epsilon = 0$      &40.3  &3.3   &47.5   &6.9   &64.2   &10.5 \\
&$\epsilon = 500$     &39.5  &3.6   &44.9   &5.0   &54.7   &7.9\\ 
&$\epsilon = 1000$    &39.3  &3.5   &43.0   &4.6   &48.0   &7.0\\ 
&$\epsilon = 1500$    &39.2  &3.5   &42.4   &3.5   &53.0   &4.8\\
&$\epsilon = 2000$    &39.3  &3.5   &40.0   &3.5   &39.8   &3.4\\
&$\epsilon = 2500$    &39.2  &3.5   &39.5   &3.4   &35.7   &3.5\\  \hline
\end{NiceTabular}
\end{center}
\end{table}

\section{Discussion of the Results}\label{discussion}
The results from Table \ref{results} clearly indicate that different risk tolerances do not affect the path-following performance and safety of the RMPC, because any worst-case collision risk would exceed the given risk tolerances. The behavior of the RMPC is only significantly influenced by the level of uncertainty. It shows that under medium and high uncertainty, the RMPC keeps an unnecessarily large distance to the object and does so at the cost of path-following performance. Contrarily, the SMPC's behavior is influenced by the risk tolerance \textit{and} the level of uncertainty. While under low uncertainty, the SMPC is not altering its behavior based on the risk tolerance, it does change its behavior for different risk tolerances under more significant uncertainties. Generally for the SMPC, with increasing risk tolerance the minimum distance to the object $d_{min}$ decreases and the path-following performance improves. Overall, the SMPC deviates less of the references and drives closer to the object than the RMPC, while not colliding. The RMPC generally appears to be strongly conservative in comparison to the SMPC.

\section{Conclusions and Future Work}\label{conclusions}

This paper provides a general stochastic definition of risk (Section \ref{riskdef}) and proposes a path-following problem with the inclusion of an explicit risk constraint (Section \ref{problem}) in the scope of autonomous driving. We proceed to derive general RMPC and SMPC formulations to solve the path-following problem (Sections \ref{rmpc} and \ref{smpc}). In the remainder of the paper, we apply both methods to a representative example, in which we simulate and compare the results of both methods in terms of safety and path-following performance (Sections \ref{implementation} - \ref{discussion}). The overarching conclusion is that the presented RMPC approach is disadvantageous for using uncertain risk measures as a constraint. The motion plans of the RMPC approach are overly conservative since the motion plan must always account for the worst-case scenario, regardless of the likelihood, and thus does not tolerate any risk. Especially with large bounds on the uncertainty, the RMPC's behavior becomes unfavorable. Contrarily, the SMPC can account for the underlying uncertainty and risk tolerances. It generates more aggressive behavior under low uncertainty and more conservative behavior under higher uncertainty, which is generally desirable.  However, the drawback of SMPC is that there are currently no proofs of the recursive feasibility of SMPC path-following schemes. If the recursive feasibility of such schemes can be proven under reasonable assumptions, safety can be genuinely guaranteed through constraint satisfaction. Our future work will focus on proving the recursive feasibility of SMPC path-following schemes. Finally, a risk term within the cost function, higher fidelity models, and increasingly complex scenarios should be investigated.

\bibliography{risk_ass}         
\bibliographystyle{ieeetr}                                                    

\appendix \label{appendix}
\section{Simulation Parameters}
\textit{1) Path:} goal point: $\boldsymbol{y}_g = (65, 5, 0)^T$, initial point: $\lambda_0 = -95$, curvature: $\kappa = 0.003$. \textit{2) Ego Vehicle:} initial configuration: $\boldsymbol{y}_{e,0} = (-10,10, 0)^T$, radius: $r_e = 1.5m$, mass: $m_e = 1000kg$, reference velocity: $v_{ref} = 3m/s$, initial input $\boldsymbol{u}_{1, e} = (0, 0)^T$. \textit{3) Object:} initial configuration: $\boldsymbol{y}_{o,0} = (5,-5,\frac{\pi}{2})^T$, radius: $r_o = 1.5m$, mass: $m_o = 1000kg$, constant inputs $\boldsymbol{u}_o = (3, 0.0001)^T$, initial standard deviation $\Sigma_0 = 0 \cdot \mathbb{I}_3$, lower bound velocity: $v_{min} = -5m/s$, upper bound velocity: $v_{max} = 5m/s$. \textit{4) Other:} sample time $T = 0.5s$, prediction horizon: $N = 6$, weighting $W = \mathbb{I}_4$,  MCS samples $J = 500$, RMPC samples $L = 40$.

\maketitle
\thispagestyle{empty}
\pagestyle{empty}

\end{document}